\documentclass[a4paper,12pt]{article}

\usepackage[T1]{fontenc}

\usepackage[latin1]{inputenc}							
\usepackage{amsmath}									

\usepackage{xcolor}
\definecolor{bl}{rgb}{0.0,0.0,0.0}
\usepackage{sectsty}
\usepackage[compact]{titlesec}
\allsectionsfont{\color{bl}\scshape\selectfont}

\makeatletter							
\def\printtitle{
    {\color{bl} \centering \huge \sc \textbf{\@title}\par}}		
\makeatother							

\title{Consistency of quasi-maximum likelihood for processes with asymetric laplacian innovation \\
		\large \vspace*{0pt} \vspace*{0pt}}

\makeatletter							
\def\printauthor{
    {\centering \small \@author}}				
\makeatother							
\author{%
	Y.BOULAROUK $^1$ \and {K.DJABALLAH} $^2$ \\
	\vspace{-10pt}
	}

\usepackage[runin]{abstract}			
\abslabeldelim{\quad}						%
\usepackage{amsmath}
\usepackage{amssymb}

\def\BibTeX{{\rm B\kern-.05em{\sc i\kern-.025em b}\kern-.08em
    T\kern-.1667em\lower.7ex\hbox{E}\kern-.125emX}}

\numberwithin{equation}{section}
\newtheorem{rmk}{Remark}[section]
\allowdisplaybreaks 
\def\1{\mbox{1\hspace{-.35em}1}} 
\def\R{\mathbb{R}}

\def\E{\mathbb{E}}

\def\R{\mathbb{R}}

\def\Z{\mathbb{Z}}

\def\v{\mbox{Var\,}}

\newtheorem{Theorem}{Theorem}

\newtheorem{Proposition}{Proposition}

\newtheorem{Lemma}{Lemma}

\begin{document}
\maketitle
\printauthor
\paragraph\
{\bf Abstract}:Strong consistency  of the quasi-maximum likelihood estimator are given for a general class of multidimensional causal processes based on asymetric laplacian innovation.
\newpage
\section{Introduction}
Since 1970 the statistical modeling changed direction, it just that the statistician comunity no longer uses structural models of Keynesian inspiration, rather stochastic models which have subsequently found a wide application in different stages of disciplines. But despite the usefulness of these processes do not cover all the phenomena and they are all(ARMA 1970,VAR 1980, ARCH 1982,....) built under the hypothesis of normality is not definitely checked and which directly affects the Likelihood function used for parameter estimation.
We give in this paper, for the first time, asymptotic properties, namely strong consistency and asymptotic normality (respectively, SC and AN for short), of the
QMLE for many multivariate models with Asymmetric Laplace errors. To establish results in a unified way, we consider almost everywhere (a.e.) solutions $X = (X_t , t \in Z)$ of equations of the type
\begin{equation}\label{eq::sys}
 X_t =M_{\theta_0}(X_{t-1},X_{t-2}, . . .)  \zeta_t +  f_{\theta_0}(X_{t-1},X_{t-2}, . . .)   t \in Z.
\end{equation}
Here, $\theta_0$   is the parameter of interest, $ M_{\theta_0}(X_{t-1},X_{t-2}, . . .)$ is a $(s\times p)$
random matrix having a.e. full rank s, $f_{\theta_0}(X_{t-1},X_{t-2}, . . 121.)$ is a $R^s$ random vector,
the ${R}^p$ random vectors $\zeta = (\zeta ^{(k)}_t )_{1\leq k\leq p}$ are independent and identically distributed satisfying standard assumptions $E[\zeta ^{(k)}_0 \zeta ^{(k')}_0 ]=0$ for $kk'$ and $E[\zeta ^{(k)^2}_0  ]=Var(\zeta ^{(k)}_0  )=1$.\\
In this study we suppose that $\zeta = (\zeta ^{(k)}_t )_{1\leq k\leq p}$ are are distributed according to an Asymetric  Laplacian law. Hence,  it has the density function $g$ given by $g(\mathbf{\zeta})=\frac{2e^{\mathbf{\zeta}^{^{\prime }}\Sigma ^{-1}m}}{%
(2\pi )^{p/2}\left\vert \Sigma \right\vert ^{1/2}}\left( \frac{\mathbf{\zeta}%
^{^{\prime }}\Sigma ^{-1}\mathbf{\zeta}}{2+m^{^{\prime }}\Sigma ^{-1}m}\right)
^{v/2}K_{v}\left( \sqrt{(2+m^{^{\prime }}\Sigma ^{-1}m)(\mathbf{\zeta}^{^{\prime
}}\Sigma ^{-1}\mathbf{\zeta})}\right)$ with m = 0 and the standard conditions defined above are checked, this function becomes :
\begin{equation}
g(\mathbf{\zeta})=\frac{2}{%
(2\pi )^{p/2} }\left( \frac{\mathbf{\zeta}%
^{^{\prime }}\mathbf{\zeta}}{2}\right)
^{v/2}K_{v}\left( \sqrt{2\mathbf{\zeta}^{^{\prime
}}\mathbf{\zeta}}\right)
\end{equation}
where $v=1-\frac{p}{2}$ and $K_v(u)$ is the modified Bessel function of the third kind given by $k_v(u)=\frac{1}{2}(\frac{u}{2})^v\int_0^\infty s^{-v-1}e^{-s-\frac{u^2}{4s}}ds, u\succ 0$ \\
Since Our sample is made up of the first n terms of an IID sequence of Laplacian random variables. The probability  density function of the vector $\zeta = (\zeta _i )_{1\leq i\leq n}$ is :
\begin{equation*}
f(\mathbf{\zeta_1,\zeta_2,....,\zeta_n})=\frac{2^n}{%
(2\pi )^{np/2} } \Pi_{t=1}^n\left[ \left( \frac{\mathbf{\zeta_t}%
^{^{\prime }}\mathbf{\zeta_t}}{2}\right)
^{v/2}K_{v}\left( \sqrt{2\mathbf{\zeta_t}^{^{\prime
}}\mathbf{\zeta_t}}\right)\right]
\end{equation*}
Through a change of variable $ X_t =M  \zeta_t + f$, we find the probability  density function of $X$ given by
$g(\mathbf{\zeta})=\frac{2 K_{v}\left( \sqrt{2\mathbf( X_t-f^t_{\theta})'(H^t)^{-1}( X_t-f^t_{\theta})}\right)}{%
(2\pi )^{np/2}det(H_t^{1/2}) }\left( \frac{\mathbf( X_t-f^t_{\theta})'(H^t)^{-1}( X_t-f^t_{\theta})}{2}\right)
^{v/2}$ and the log likelihood function is
\begin{equation}
L_n(\Theta)=\sum_{t=1}^n  q_t(\Theta)
\end{equation}

With
\begin{align}
	\begin{split}
	q_t(\Theta) 	&=log\left(K_v \sqrt{ 2 ( X_t-f^t_{\Theta})'(H^t)^{-1}( X_t-f^t_{\Theta})}\right)  + \frac{v}{2} log\left({  ( X_t-f^t_{\Theta})'(H^t)^{-1}( X_t-f^t_{\Theta})}\right)\\
					&-\frac{1}{2}log \left( det (H^t_\Theta)\right)\\
	\end{split}					
\end{align}

\begin{equation}\label{q::mle}
\widehat \theta_n:=\underset{\theta \in \Theta}{\mbox{Argmax }}\widehat
L_n(\theta).
\end{equation}

\subsection{Definition of the parameter sets $\Theta(r)$ and $\widetilde \Theta(r)$}
In proposition \ref{stationarity} below we provide the existence of a stationary solution of the general model \eqref{eq::sys}. Two conditions of different types are used: the first one is a Lipschitz condition on the functions $f$ and $M$ in \eqref{eq::sys}, the second one is a restriction on the set of the parameters.\\

Let us assume that for any $\theta\in\R^d$, $x\mapsto f_\theta(x)$ and $x\mapsto M_\theta (x)$ are Borel functions on $({\R^m}) ^\infty$ and that $\mbox{Rank }M_\theta (x)= m$ for all $x\in({\R^m}) ^\infty$. Assume that there exist two sequences $(\alpha_j(f,\theta))_{j\ge1}$ and $(\alpha_j(M,\theta))_{j\ge1}$ satisfying, for all $x$,
$y$ in $({\R^m}) ^{\infty}$,
$$\left \{
\begin{array}{lll}
\|f_{\theta}(x)-f_{\theta}(y)\|
&\le &
\sum_{j=1}^\infty \alpha_j(f,\theta)\|x_j-y_j\|,\\
\|M_{\theta}(x)-M_{\theta}(y)\|
&\le &\sum_{j=1}^\infty \alpha_j(M,\theta)\|x_j-y_j\|.
\end{array}
\right .
$$
We can define the set
\begin{equation}
\label{propINF}\Theta(r)=\left\{\theta\in\R^d~\Big/~ \sum_{j=1}^\infty
\alpha_j(f,\theta)+\left(\E\|\xi_0\|^r\right)^{1/r}\sum_{j=1}^{\infty} \alpha_j(M,\theta)<1\right\}.
\end{equation}
This set depends on the distribution of $\xi_0$ via the moments $\E\|\xi_0\|^r$. But thanks to the fact that $\E\big [\xi_0^{(k)}\xi_0^{(k')}\big ]=0$ for $k\neq k'$ and $\E
\big [{\xi_0^{(k)}}^2\big ]=\v(\xi_0^{(k)})=1$ the set $\Theta(2)$ simplifies:
$$
\Theta(2)=\left\{\theta\in\R^d~\Big/~  \sum_{j=1}^\infty
\alpha_j(f,\theta)+\sqrt{p}\sum_{j=1}^{\infty} \alpha_j(M,\theta)<1\right\}.
$$
\begin{Proposition}\label{stationarity}
If $\theta_0\in \Theta(r)$ for some $r\ge1$ there exists a unique
causal ($X_t$ is independent  of $(\xi_i)_{i> t}$ for $t\in\Z$) solution $X$ to the equation \eqref{eq::sys} which is
stationary and ergodic  and satisfies
$\E\big \|X_0\big\|^r <\infty$.
\end{Proposition}
\subsection{Uniform assumptions on $\Theta$}
Fix some compact subset $\Theta$ of $\R^d$. For  any
sequences $x$, $y$ of $({\R^m}) ^{\infty}$, the functions $\theta
\mapsto f_{\theta}(x)$ and $\theta \mapsto M_\theta(x)$ are
assumed to be continuous on $\Theta$. Assume that $\|f_\theta(0)\|_{\Theta}<\infty$ and $
\|M_\theta(0)\|_{\Theta} <\infty$. To settle the
assumptions in a short way, let us introduce the generic symbol $\Psi$ for any of the functions
$f$, $M$ or $H$.
\begin{description}
\item[(A1($\Psi$))] Let $\alpha_j(\Psi)=\sup_{\theta\in\Theta}\alpha_j(\Psi,\theta)$
be such that $\sum_{j\ge1} \alpha_j(\Psi)<\infty$.
\item[(A2)] There exists $\underline H>0$ such that
$\inf_{\theta\in\Theta}\det \big(H_\theta(x) \big ) \ge
\underline H$ for all $x\in ({\R^m}) ^{\infty}$.
\item[(A3($\Psi$))] The function $\theta \in \Theta  \mapsto
\Psi_\theta(x)$ is $2$ times continuously differentiable for all
$x\in(\R^m)^{\infty}$ and
$$ \Big \|\frac {\partial
\Psi_\theta(0)}{\partial \theta} \Big \|_\Theta+\Big \|\frac
{\partial^2 \Psi_\theta(0)}{\partial \theta\partial \theta'} \Big
\|_\Theta< \infty.$$
Moreover assume that there exist two integrable sequences
$\big(\alpha^{(i)}_j(\Psi)\big)_{j\ge1}$, $i=1,2$, such that for all $x$, $y\in(\R^m)^{\infty}$
\begin{eqnarray*}\Big \|\frac
{\partial \Psi_\theta(x)}{\partial \theta
 } -\frac {\partial \Psi_\theta(y)}{\partial \theta} \Big \|_\Theta &\le& \displaystyle  \sum_{j=1}^\infty
\alpha^{(1)}_j(\Psi)
\|x_j-y_j\|,\\
\Big \|\frac
{\partial^2 \Psi_\theta(x)}{\partial \theta\partial \theta'
 } -\frac {\partial^2 \Psi_\theta(y)}{\partial \theta\partial \theta'} \Big \|_\Theta &\le& \displaystyle  \sum_{j=1}^\infty
\alpha^{(2)}_j(\Psi)
\|x_j-y_j\|.
\end{eqnarray*}
If $\Psi=H$, $\|x_j-y_j\|$ in the RHS terms is replaced with $\|x_jx_j'-y_jy_j'\|$.
\end{description}
The last assumption on the derivatives is just needed for the asymptotic normality of the QMLE.
\subsection{Identifiability and variance conditions}
We assume the same identifiability condition as in Jeantheau \cite{jen}:
\begin{description}
\item[(Id)] For all $\theta\in \Theta$,
($f^t_\theta=f^t_{\theta_0}$ and
$H^t_\theta=H^t_{\theta_0}$ a.s.) $\Rightarrow\theta=\theta_0$.\\
\item[(Var)] One of the families $({\partial f_{\theta_0}^t}/{\partial
\theta_i})_{1\le i\le d}$ or $({\partial
H_{\theta_0}^t}/{\partial \theta_i})_{1\le i\le d}$ is a.e.
linearly independent, where:
$$
\frac{\partial f_\theta^t}{\partial \theta}:=
\frac{\partial f_\theta}{\partial \theta}(X_{t-1},\ldots)\mbox{ and }
\frac{\partial H_\theta^t}{\partial \theta}:=
\frac{\partial H_\theta}{\partial \theta}(X_{t-1},\ldots).
$$
\end{description}
The condition {\bf (Var)} is needed for ensuring finiteness of the asymptotic variance in the result on asymptotic normality.
\begin{Proposition}
Let be $k_v(u)$the bessel function of third kind. For all real x,y there exist A,B constants which satisfies :
\begin{itemize}
   \item $(ii)   \*sup_{\theta \in \Theta}|k_v(x)|<A u^v$
   \item $(ii)   \sup_{\theta \in \Theta}|k_v(y)-k_v(x)|<B|y-x|      \  \ \        ,\forall v\leq 0$
   \item $(iii)  \sup_{\theta \in \Theta}|log[k_{\frac{1}{2}}(y)]-log[k_{\frac{1}{2}}(x)|]<B|y-x|$
 \end{itemize}  
\end{Proposition}
proof :\\
(i) We have :\\
\begin{align*}
	\begin{split}
k_v(u)=\frac{1}{2} (\frac{u}{2})^v\int_0^\infty t^{-v-1}e^{-t-\frac{u^2}{4t}}dt	\end{split}					
\end{align*}
\begin{align*}
	\begin{split}
|k_v(u)|	&=\frac{1}{2} |\frac{u}{2}|^v\int_0^\infty t^{-v-1}e^{-t-\frac{u^2}{4t}}dt   \\
				&   =\frac{1}{2^{v+1}}|u|^v\int_0^\infty t^{-v-1}e^{-t-\frac{u^2}{4t}}dt\\
				& =\frac{1}{2^{v+1}}|u|^v\int_0^\infty t^{-v-1}e^{-t}e^{-\frac{u^2}{4t}}dt  \\
				&\leq \frac{1}{2^{v+1}}|u|^v\int_0^\infty t^{-v-1}e^{-t}dt\leq A u^v\ \ \forall v<0   \\
	\end{split}					
\end{align*}
(ii)
\begin{align*}
	\begin{split}
k_v(u)	&=\frac{1}{2}(\frac{u}{2})^v\int_0^\infty t^{-v-1}e^{-t-\frac{u^2}{4t}}dt\\
				& =\frac{1}{2^{v+1}}\int_0^\infty t^{-v-1}e^{-t}{u}^v e^{-\frac{u^2}{4t}}dt     \\
	\end{split}					
\end{align*}

Let put $g(u)={u}^v e^{-\frac{u^2}{4t}}$  ($g'(u)={u}^{v-1} e^{-\frac{u^2}{4t}}(v-\frac{u^2}{4t})$
$\lim_{x\rightarrow 0}g'(u)=a < 0$ and $\lim_{u\rightarrow \infty } g'(u)=0$) \\
$g''(u)={u}^{v} e^{-\frac{u^2}{4t}}(v-\frac{u^2}{2t}+1)$
$g''(u)=0 \Rightarrow u=\pm \sqrt{- 2t(v+1)}$\\ *
a simple study of $g'$ shows that it is bounded and so $g$ is lipchetzian
which implies that $|g(y)-g(x)| \leq c |y-x|$\\

\begin{align*}
	\begin{split}
|k_v(y)-k_v(x)|    &\leq \frac{1}{2^{v+1}}\int_0^\infty t^{-v-1}e^{-t}|{y}^v e^{-\frac{y^2}{4t}}-{x}^v e^{-\frac{x^2}{4t}}|dt + \\
   		 &\leq  \frac{1}{2^{v+1}}\int_0^\infty t^{-v-1}e^{-t}c |y-x|dt\\
		    &=\frac{c}{2^{v+1}}|y-x|\int_0^\infty t^{-v-1}e^{-t}dt
	\end{split}					
\end{align*}
or :\\
\begin{align*}
	\begin{split}
\int_0^\infty t^{-v-1}e^{-t}dt
=\int_0^1 t^{-v-1}e^{-t}dt + \int_1^\infty t^{-v-1}e^{-t}dt<\infty \  \ since  \ v<0
	\end{split}					
\end{align*}
\begin{align*}
	\begin{split}
|k_0(y)-k_0(x)|    &\leq \int_0^1\frac{|cos(yt)-cos(xt)|}{\sqrt {t^2+1}}dt +\frac{1}{2^{v+1}}\int_1^\infty t^{-v-1}e^{-t}|{y}^v e^{-\frac{y^2}{4t}}-{x}^v e^{-\frac{x^2}{4t}}|dt  \\
   		 &\leq  \int_0^1\frac{c |y-x|t}{\sqrt {t^2+1}}dt+\frac{1}{2^{v+1}}\int_1^\infty t^{-v-1}e^{-t}c' |y-x|dt\\
		    &=c |y-x|\left[ \int_0^1\frac{t}{\sqrt {t^2+1}}dt+\frac{1}{2^{v+1}}\int_1^\infty t^{-v-1}e^{-t}dt\right]\\
		 &=c |y-x|\left[\sqrt{2}+\frac{1}{2^{v+1}}\int_1^\infty t^{-v-1}e^{-t}dt\right]<\infty
	\end{split}					
\end{align*}
so
\begin{align*}
	\begin{split}
\sup_{\theta \in \Theta}|k_v(y)-k_v(x)|\leq B|y-x|    \      \   \ \forall v\leq 0
	\end{split}					
\end{align*}
\\
(iii) by definitin we have
$k_v(u)=\frac{1}{2}(\frac{u}{2})^v\int_0^\infty t^{-v-1}e^{-t-\frac{u^2}{4t}}dt, u\succ 0$
in particular
\begin{align*}
	\begin{split}
k_\frac{1}{2}(u)=	&\sqrt{\frac{\pi}{2u}  }e^{-u}\\
		& log[k_\frac{1}{2}(u)]=\frac{1}{2}[log\pi - log(2u)]-u
	\end{split}					
\end{align*}
else
\begin{align*}
	\begin{split}
| log[k_\frac{1}{2}(y)- log[k_\frac{1}{2}(x)]|		&=|(\frac{1}{2}[log\pi - log(2y)]-y)-( \frac{1}{2}[log\pi - log(2x)]-x)|\\
	&=\frac{1}{2}| log(2x)- log(2y)  + 2(x-y)|\leq c |y-x|
	\end{split}					
\end{align*}
\begin{rmk}
   This last proposition implies that :$\left| log[k_v(y)- log[k_v(x)]\right|\leq c |y-x|$ for any possible values of $v$.
\end{rmk}
\begin{Theorem}
Assume that $\theta_0 \in \Theta(2)\cap \Theta$ and let $X$ be the stationary solution
of \eqref{eq::sys}. If  $\theta_0 \in  \Theta$ , a compact set of $R^d$ such that assumptions $(D(\Theta)), (Id(\Theta)),
(A0(f,\Theta)) and (A0(M,\Theta)) [or (A0(H,\Theta))]$ hold with :
\begin{equation}
\alpha_j^0(f,\Theta)+\alpha_j^0(M,\Theta)+\alpha_j^0(H,\Theta)=\mathcal{O}(j^-l) for some  \   l\succ {2}
\end{equation}
then the QMLE  $\widehat{\theta_n}$ defined in \eqref{q::mle} is SC; that is,  $\widehat{\theta_n}   \rightarrow        \theta_0$ a.s.
\end{Theorem}
\begin{Lemma}
 Assume that $\theta_0 \in \Theta(r)$ for $r\geq 2$ and $X$ is the stationary solution
of (1). Let $\Theta$ be a compact set of $R^d$ :\\
\begin{enumerate}
  \item If $(A_0(f,\Theta))$ holds, then $\forall\theta \in \Theta, f_\theta^t \in L^r(C{(\Theta,\mathcal{R}^m))}$ and
\begin{equation*}
E[\|\widehat{f}_\theta^t-{f}_\theta^t\|_\Theta^r]\leq E\|X_0\|^r \left(\|\alpha_j(f)\|_\Theta\right)^r \ for \ all   \  t\in N^*;
\end{equation*}
  \item If $(A_0(M,\Theta))$ holds, then $\forall\theta \in \Theta, H_\theta^t \in L^{r/2}(C{(\Theta,\mathcal{M}_m))}$ and there exists $C>0$ not depandig on t such that
\begin{equation*}
E[\|\widehat{H}_\theta^t-{H}_\theta^t\|_\Theta^{r/2}]\leq C \left(\|\alpha_j^{(0)}(M,\Theta)\|_\Theta\right)^{r/2}  \ for \ all   \  t\in N^*;
\end{equation*}
  \item If $(A_0(H,\Theta))$ holds, then $\forall\theta \in \Theta, H_\theta^t \in L^{r/2}(C{(\Theta,\mathcal{M}_m))}$ and
\begin{equation*}
E[\|\widehat{H}_\theta^t-{H}_\theta^t\|_\Theta^{r/2}]\leq E\|X_0\|^r \left(\|\alpha_j^{(0)}(H,\Theta)\|_\Theta\right)^{r/2}  \ for \ all   \  t\in N^*;
\end{equation*}
\end{enumerate}
Moreover, under any of the two last conditions and with $(D(\Theta))$, $H_\theta^t$ is an invertible
matrix and $\|(\widehat{H}_\theta^t)^{-1}\|_\Theta\leq\underline{H}^{-1/m}$. 
\end{Lemma}
$\mathbf{Proof\  of \ Lemma 1.}$  See Lemma1 of Bardet and Weitenberger (\cite{bardet}).\\
$\mathbf{Proof\  of \ Theorem 1.}$
The proof of the theorem is divided into two parts.
In (i), a uniform (in $\theta$) law of large numbers on $(\widehat{q}_t )_{t\in N^*}$ [defined in (1.3)] is established.
In (ii), it is proved that $L(\theta) := E(q_0 (\theta))$ has a unique maximum in $\theta_0$.
Those two conditions lead to the consistency of $\widehat{\theta}_n$.
(i) Using \textit{Proposition 1}, with $q_t = G(X_t,X_{t-1}, . . .)$, one deduces that $(q_t )_{t\in Z}$
[defined in (3)] is a stationary ergodic sequence. From Straumann and Mikosch
\cite{STRAU}, we know that, if $(v_t )_{t\in Z}$ is a stationary ergodic sequence of random elements
with values in $C(\Theta,R^m)$, then the uniform (in $\theta\in \Theta$) law of large numbers is implied
by$ E||v_0||_\Theta\leq \infty$.As a consequence,  $X = (q_t , t \in Z)$ satisfies a uniform $ (in \theta \in \Theta)$ strong law of large numbers as soon as
$ E[sup_{\theta}|q_t (\theta)|] \prec \infty$. But, from the inequality $ log(x)\leq x$, for all $x \in ]0,\infty[$ and Lemma 1, for all $t \in Z$,\\
\begin{align*}
	\begin{split}
|q_t(\theta)|	&\leq \left|log \left[ K_v\left( \sqrt{ 2 ( X_t-f^t_{\theta})'(H^t)^{-1}( X_t-f^t_{\theta})}\right)\right] \right| \\
				&+\left| \frac{v}{2}{  ( X_t-f^t_{\theta})'(H^t)^{-1}( X_t-f^t_{\theta})}-\frac{1}{2}log \left( det (H^t_\theta)\right)\right|
	\end{split}					
\end{align*}
Note :  $S(\Theta)= \frac{v}{2}    {  ( X_t-f^t_{\Theta})'(H^t)^{-1}( X_t-f^t_{\Theta})}-\frac{1}{2}log \left( det (H^t_\Theta)\right)$
\begin{align*}
	\begin{split}
\left|S(\Theta)\right| \leq \frac{v}{2}\frac{|| X_t-f^t(\theta)||^2}{(\underline{H})^\frac{1}{m}}+\frac{m}{2}\left|\frac{1}{m}log(\underline{H})+\frac{|| H^t_(\theta)||}{(\underline{M})^\frac{1}{m}}-1\right|\\
	\end{split}					
\end{align*}
	
\begin{equation}\label{sd:ss}
\Longrightarrow
sup_{\theta \in \Theta}   \left|S(\theta)\right| \leq \frac{v}{2}\frac{|| X_t-f^t(\theta)||^2_\Theta}{\underline{H}^\frac{1}{m}}+\frac{1}{2}|log\underline{H}|+\frac{m}{2}\frac{|| H^t_\theta||_\Theta}{\underline{H}^\frac{1}{m}}
\end{equation}
But, $\forall t \in Z,E||X_t||\prec \infty$ (see Proposition 1) and $E||f^t(\theta)||^r_\Theta+E||H^t(\theta)||^r_\Theta\prec \infty$
(see Lemma 1). As a consequence, the right-hand side of \eqref{sd:ss} has a finite first moment
.Therefore to proof that $q_t$ have a finite first order moment, we have to proof that $E \left(\sup_{\theta \in \Theta}\left|log\left[K_v \sqrt{ 2 ( X_t-f^t_{\Theta})'(H^t)^{-1}( X_t-f^t_{\Theta})}\right)\right]\right|  \prec \infty$.\\
\\
By the result of \textit{Proposition 2} we  have $\sup_{\theta \in \Theta}|k_v(u)|<A u^v$\\
\begin{align*}
	\begin{split}
\left|log\left[K_v \sqrt{ 2 ( X_t-f^t_{\theta})'(H^t)^{-1}( X_t-f^t_{\theta})}\right]\right|  	&\leq \left|K_v \sqrt{ 2 ( X_t-f^t_{\theta})'(H^t)^{-1}( X_t-f^t_{\theta})}\right|\\
	&\leq A'\left[{ 2 ( X_t-f^t_{\theta})'(H^t)^{-1}( X_t-f^t_{\theta})}\right]^{\frac{v}{2}}\\
	&\leq A' \frac{|| X_t-f^t(\theta)||_\Theta}{\underline{H}^\frac{1}{2 m}}
	\end{split}					
\end{align*}
\begin{equation*}
\Longrightarrow
\sup_{\theta \in \Theta}\left|log \left[ K_v\left( \sqrt{ 2 ( X_t-f^t_{\Theta})'(H^t)^{-1}( X_t-f^t_{\Theta})}\right)\right]\right| 	 \leq A'  \frac{|| X_t-f^t(\theta)||_\Theta}{\underline{H}^\frac{1}{2 m}}
\end{equation*}
whence :
\begin{align*}
	\begin{split}
E\sup_{\theta \in \Theta}\left|log\left[K_v \sqrt{ 2 ( X_t-f^t_{\Theta})'(H^t)^{-1}( X_t-f^t_{\Theta})}\right]\right|  	 < \infty
	\end{split}					
\end{align*}
and, therefore,
\begin{equation*}
E \left[sup_{\theta \in \Theta}\left| q_t(\theta) \right|\right]\prec\infty
\end{equation*}
The uniform strong law of large numbers for $(q_t (\theta))$ follows; hence,
\begin{equation*}
||\frac{L_n(\theta)}{n}-L(\theta)||_\Theta\rightarrow 0 \  a.s  \   \   \   \   with   \  L(\theta):=E[q_0(\theta)]
\end{equation*}
Now, one shows that $\frac{1}{n}||\widehat{L_n}-L_n||_\Theta\rightarrow 0  \  a.s$. Indeed, for all $\theta \in \Theta$ and $t\in N^*$,
Let put
\begin{align*}
{A_t}(\theta)=\frac{1}{2}[v {  ( X_t-f^t_{\theta})'(H^t)^{-1}( X_t-f^t_{\theta})}+det(H^t_{\theta})]\\
{B_t}(\theta)= log\left(K_v \sqrt{ 2 ( X_t-f^t_{\theta})'(H^t)^{-1}( X_t-f^t_{\theta})}\right)
\end{align*}
\begin{align}
	\begin{split}
|\widehat{q_t}(\theta)-q_t(\theta)|	&\leq |\widehat{A_t}(\theta)-A_t(\theta)|+|\widehat{B_t}(\theta)-B_t(\theta)|
	\end{split}					
\end{align}
\begin{align}
	\begin{split}
|\widehat{A_t}(\theta)-A_t(\theta)|	&\leq \frac{v}{2} \left| ( X_t-\widehat f^t_{\theta})'(\widehat H^t)^{-1}( X_t-\widehat f^t_{\theta})- ( X_t-f^t_{\theta})'(H^t)^{-1}( X_t-f^t_{\theta})\right|  \\
				&+\frac{1}{2}\left| det (\widehat H^t_\theta)- det ( H^t_\theta)\right|\\
				&\leq\frac{1}{2|C|}\left| det (\widehat H^t_\theta)- det ( H^t_\theta)\right| +\frac{v}{2}  ( X_t-\widehat f^t_{\theta})'[(\widehat H^t)^{-1}-(H^t)^{-1}]( X_t-\widehat f^t_{\theta}) \\
				&+ \frac{v}{2}(2 X_t-\widehat f^t_{\theta}-f^t_{\theta})'(H^t)^{-1}(f^t_{\theta}-\widehat f^t_{\theta})\\
				&\leq\frac{1}{2}H^{-1}\left|\left| det (\widehat H^t_\theta)- det ( H^t_\theta)\right|\right|_{\Theta}+\frac{v}{2} \left( \left|\left| X_t\right|\right|+\left|\left|\widehat f^t_{\theta}\right|\right|_{\Theta})\left|\left|(\widehat H^t)^{-1}-(H^t)^{-1}\right|\right|_\Theta\right)  \\
				&+ \frac{v}{2}\left((2 \left|\left|X_t\right|\right|+||\widehat f^t_{\theta}||_\Theta+\left|\left|f^t_{\theta}\right|\right|_\Theta)       \left|\left|(H^t)^{-1}\right|\right|_{\Theta}            \left|\left|f^t_{\theta}-\widehat f^t_{\theta}\right|\right|_{\Theta}\right)
	\end{split}					
\end{align}
on the one hand we have,
$\left|\left|(\widehat H_\theta^t)^{-1}-(H_\theta^t)^{-1}\right|\right|_\Theta \leq \left|\left|(\widehat H_\theta^t)^{-1}\right|\right|_\Theta  \left|\left|\widehat H_\theta^t-H_\theta^t\right|\right| \left|\left|( H^t)^{-1}\right|\right|_\Theta $
 on the other hand, for invertible matrix $A\in M_m(R)$ , and $H\in  M_m(R)$,
\begin{align*}
	\begin{split}
det(H_\theta^t) = det(\widehat H_\theta^t) +det(\widehat H_\theta^t)  . Tr\left((\widehat H_\theta^t)^{-1})'||(\widehat H_\theta^t)^{-1}-(H_\theta^t)^{-1}||\right)+o(||(\widehat H^t)^{-1}-(H^t)^{-1}||),
	\end{split}					
\end{align*}
where  $\left|Tr\left(((\widehat H_\theta^t)^{-1})'||(\widehat H_\theta^t)^{-1}-(H_\theta^t)^{-1}||\right)\right|\leq \left|\left|(\widehat H^t)^{-1}\right|\right|_\Theta \left|\left|\widehat H_\theta^t-H_\theta^t\right|\right| $. Using the relation $\|(H_\theta)^{-1}\|\geq\underline{H}^{-m} $ for all $t\in Z$, there exists $C>0$ not depending on $t$ , such that inequality (8) becomes
\begin{align*}
	\begin{split}
\sup_{\theta \in \Theta}|\widehat{A_t}(\theta)-A_t(\theta)|	\leq C\left( \left|\left|X_t \right|\right|+\left|\left|\widehat{f}_\theta^t \right|\right|_{\Theta}+\left|\left|{f}_\theta^t \right|\right|_{\Theta}\right)\left(\|\widehat{H}_\theta^t-{H}_\theta^t\|_\Theta+\|\widehat{f}_\theta^t-{f}_\theta^t\|_\Theta\right)
	\end{split}					
\end{align*}
\begin{align*}
	\begin{split}
|\widehat{B_t}(\theta)-B_t(\theta)|	&\leq  \left|K_v \sqrt{ 2 ( X_t-f^t_{\widehat{\Theta}})'(\widehat{H}^t)^{-1}( X_t-f^t_{\widehat{\Theta}})}-K_v \sqrt{ 2 ( X_t-f^t_{\Theta})'(H^t)^{-1}( X_t-f^t_{\Theta})}\right|\\
				&  \leq A \left| \sqrt{ 2 ( X_t-f^t_{\widehat{\Theta}})'(\widehat{H}^t)^{-1}( X_t-f^t_{\widehat{\Theta}})}- \sqrt{ 2 ( X_t-f^t_{\Theta})'(H^t)^{-1}( X_t-f^t_{\Theta})}\right|\\
				&\leq \sqrt{2} A  \sqrt{  ( X_t-f^t_{\widehat{\Theta}})'(\widehat{H}^t)^{-1}( X_t-f^t_{\widehat{\Theta}})-  ( X_t-f^t_{\Theta})'(H^t)^{-1}( X_t-f^t_{\Theta})} \\
				&\leq \sqrt{2}  A \left({ ( X_t-\widehat f^t_{\theta})'[(\widehat H^t)^{-1}-(H^t)^{-1}]( X_t-\widehat f^t_{\theta})+ (2 X_t-\widehat f^t_{\theta}-f^t_{\theta})'(H^t)^{-1}(f^t_{\theta}-\widehat f^t_{\theta})})\right)^\frac{1}{2}		
	\end{split}					
\end{align*}
Following the same approach for A found :
\begin{align*}
	\begin{split}
\sup_{\theta \in \Theta}|\widehat{B_t}(\theta)-B_t(\theta)|	\leq C' \left( \left|\left|X_t \right|\right|+\left|\left|\widehat{f}_\theta^t \right|\right|_{\Theta}+\left|\left|{f}_\theta^t \right|\right|_{\Theta}\right)^\frac{1}{2}\left(\|\widehat{H}_\theta^t-{H}_\theta^t\|_\Theta+\|\widehat{f}_\theta^t-{f}_\theta^t\|_\Theta\right)^\frac{1}{2}
	\end{split}					
\end{align*}
From the Holder and Minkowski inequalities and by virtue of $3/2 = 1+ 1/2$,
\begin{align}
	\begin{split}
E[\sup_{\theta \in \Theta}|\widehat{A_t}(\theta)-A_t(\theta)|]	&\leq C\left(E[ \left|\left|X_t \right|\right|+\left|\left|\widehat{f}_\theta^t \right|\right|_{\Theta}+\left|\left|{f}_\theta^t \right|\right|_{\Theta}]^2\right)^{1/2}\\
			& \times\left(E[\|\widehat{H}_\theta^t-{H}_\theta^t\|_\Theta]+E[\|\widehat{f}_\theta^t-{f}_\theta^t\|_\Theta]^2\right)^{1/2}\\
			&\leq C_*\left(E[\|\widehat{H}_\theta^t-{H}_\theta^t\|_\Theta]+E[\|\widehat{f}_\theta^t-{f}_\theta^t\|_\Theta]^2\right)^{1/2}
	\end{split}					
\end{align}
\begin{align}
	\begin{split}
E[\sup_{\theta \in \Theta}|\widehat{B_t}(\theta)-B_t(\theta)|^{1/2}]	 &\leq C'\left(E \left|\left|X_t \right|\right|+\left|\left|\widehat{f}_\theta^t \right|\right|_{\Theta}+\left|\left|{f}_\theta^t \right|\right|_{\Theta}\right)^{1/2}\\
			& \times\left(E[\|\widehat{H}_\theta^t-{H}_\theta^t\|_\Theta]+E[\|\widehat{f}_\theta^t-{f}_\theta^t\|_\Theta]\right)^{1/2}\\
	&\leq C_*'\left(E[\|\widehat{H}_\theta^t-{H}_\theta^t\|_\Theta]+E[\|\widehat{f}_\theta^t-{f}_\theta^t\|_\Theta]\right)^{1/2}
	\end{split}					
\end{align}
with $C_*>0,C'_*>0 $ not depending on $\theta$ and $t$ . Now, consider, for $n \in N$,
\begin{align*}
	\begin{split}
S_n = \sum_{t=1}^n\frac{1}{t}\sup|\widehat{q_t}(\theta)-q_t(\theta)|.
	\end{split}					
\end{align*}
Applying the Kronecker lemma (see Feller \cite{Feller}, page 238), if $\lim_n\rightarrow \infty Sn < \infty$ a.s., then $\frac{1}{n}\|\widehat{L_n}-L_n\|\rightarrow 0 a.s$. Following Feller's arguments, it remains to show
that, for all $\varepsilon>0$,
\begin{align*}
	\begin{split}
P(\forall n \in {N},\exists m \ such \  that \ |S_m-S_n|>\varepsilon)={P}(A)=0.
	\end{split}					
\end{align*}
Let $\varepsilon>0$, and denote
\begin{align*}
	\begin{split}
A_{m,n} := \{|S_m - S_n| > \varepsilon\}
	\end{split}					
\end{align*}
for $m >n$. Notice that $A =\bigcap_{n\in N}\bigcup_m>n A_m,n$. For $n\in N^*$, the sequence of
sets $(A_{m,n})_m>n$ is obviously increasing, and, if $A_n=\bigcup_m>n A_{m,n}$, then
$lim_{m\rightarrow\infty}P(A_{m,n}) = P(A_n)$. Observe that $(A_n)_{n\in N}$ is a decreasing sequence of sets and, thus,
\begin{align*}
	\begin{split}
lim_{n\rightarrow\infty} lim_{m\rightarrow\infty}P(A_{m,n}) = lim_{n\rightarrow\infty} P(A_n) = P(A).	
	\end{split}					
\end{align*}
It remains to bound $P(A_{m,n})$. From the $Bienyam\acute{e} - Chebyshev$ inequality,
\begin{align*}
	\begin{split}
P(A_{m,n})	&=  P\left(\sum_{t=n+1}^m \frac{1}{t} \sup_{\theta\in \Theta} |\widehat{q_t}(\theta)-q_t(\theta)|>\varepsilon\right)\\
		&\leq\frac{1}{\varepsilon^{\frac{2}{3}}} E\left[\left(\sum_{t=1}^n \frac{1}{t} \sup_{\theta\in \Theta} |\widehat{q_t}(\theta)-q_t(\theta)|\right)^{\frac{2}{3}}\right]\\
		&\leq \frac{1}{\varepsilon^{\frac{2}{3}}} \sum_{t=n+1}^m \frac{1}{t^{\frac{2}{3}}} E\left(\sup_{\theta\in \Theta} |\widehat{q_t}(\theta)-q_t(\theta)|^{\frac{2}{3}}\right).
	\end{split}					
\end{align*}
Using (9) and condition (5), since  $l> 3/2$, there exists $C >0$ such that
\begin{equation*}
\left(\sum_{j=t}^\infty \alpha_j^0(f,\Theta)+\alpha_j^0(M,\Theta)+\alpha_j^0(H,\Theta)\right) ^{\frac{2}{3}}\leq \frac{C}{t^{2(l-1)/3}}
\end{equation*}
Thus, ${t^{\frac{-2}{3}}} E\left( \sup_{\theta\in \Theta} |\widehat{q_t}(\theta)-q_t(\theta)|^{\frac{2}{3}}\right)<C(t^{-2l/3})$ for some $C >0$, and
\begin{equation*} \sum_{t=1}^m \frac{1}{t^{\frac{2}{3}}} E\left(\sup_{\theta\in \Theta} |\widehat{q_t}(\theta)-q_t(\theta)|^{\frac{2}{3}}\right)<\infty    \    \    \    \     \ a.s\    \      l>3/2
\end{equation*}
Thus, $lim_{n\rightarrow\infty}lim_{m\rightarrow\infty}P(A_{m,n}) \rightarrow0$ and $\frac{1}{n}\|\widehat{L_n}-L_n\|\rightarrow 0 a.s$\\
(ii) En cours.......\\

\end{document}